\input ./style/arxiv-general.cfg
\documentclass[aos,MSNbibl,secthm,seceqn,dvips]{arximspdf}
\makeatletter
   \@ifpackageloaded{graphicx}{}{\usepackage{graphicx}}
\makeatother
\usepackage{dcolumn}
\usepackage{algorithm}
\usepackage{algpseudocode}

\doi{10.1214/15-AOS1348}
\volume{43}
\issue{6}
\pubyear{2015}
\firstpage{2484}
\lastpage{2506}
\docsubty{FLA}

\makeatletter
\newcommand{\rrvert}{\vert}
\newcommand{\llvert}{\vert}
\newcolumntype{d}[1]{D{.}{.}{#1}}
\newtheorem{prop}{Proposition}[section]
\newtheorem{lemma}{Lemma}[section]
\makeatother

\begin{document}
\begin{frontmatter}

\title{Coupling methods for multistage sampling}
\runtitle{Coupling methods for multistage sampling}

\begin{aug}
\author[A]{\fnms{Guillaume}~\snm{Chauvet}\corref{}\ead[label=e1]{chauvet@ensai.fr}}
\runauthor{G. Chauvet}
\affiliation{CREST-ENSAI}
\address[A]{CREST-ENSAI\\
Rue Blaise Pascal\\
Campus de Ker Lann \\
35170 Bruz \\
France \\
\printead{e1}}
\end{aug}

%
\received{\smonth{5} \syear{2015}}
%
\revised{\smonth{5} \syear{2015}}

\begin{abstract}
Multistage sampling is commonly used for household surveys when there
exists no sampling frame, or when the population is scattered over a
wide area. Multistage sampling usually introduces a complex dependence
in the selection of the final units, which makes asymptotic results
quite difficult to prove. In this work, we consider multistage sampling
with simple random without replacement sampling at the first stage, and
with an arbitrary sampling design for further stages. We consider
coupling methods to link this sampling design to sampling designs where
the primary sampling units are selected independently. We first
generalize a method introduced by [\textit{Magyar Tud. Akad. Mat. Kutat\'o Int. K\"ozl.}
\textbf{5} (1960) 361--374]
to get a coupling with
multistage sampling and Bernoulli sampling at the first stage, which
leads to a central limit theorem for the Horvitz--Thompson estimator.
We then introduce a new coupling method with multistage sampling and
simple random with replacement sampling at the first stage. When the
first-stage sampling fraction tends to zero, this method is used to
prove consistency of a with-replacement bootstrap for simple random
without replacement sampling at the first stage, and consistency of
bootstrap variance estimators for smooth functions of totals.
\end{abstract}

\end{frontmatter}

\section{Introduction} \label{secintro}

Multistage sampling is widely used for household and health surveys
when there exists no sampling frame, or when the population is
scattered over a wide area. Three or more stages of sampling may be
commonly used. For example, the third National Health and Nutrition
Survey (NHANES III) conducted in the United States involved four stages
of sampling, with the selection of counties as Primary Sampling Units
(PSUs), of segments as Secondary Sampling Units (SSUs) inside the
selected counties, of households as Tertiary Sampling Units (TSUs)
inside the selected segments, and of individuals inside the selected
households, for example, \cite{ezzhofjudmasmoo92}. A detailed
treatment of multistage sampling may be found in textbooks like \cite{coc77,sarswewre92} or \cite{ful11}.

Multistage sampling introduces a complex dependence in the selection of
the final units, which makes asymptotic properties difficult to prove.
In this work, we make use of coupling methods (see \cite{tho00}) to
link multistage sampling designs to sampling designs where the primary
sampling units are selected independently. The method basically
consists in generating a random vector $(X_t,Z_t)^{\top}$ with
appropriate marginal laws, and so that $E(X_t-Z_t)^2$ is smaller than
the rate of convergence of $X_t$. In this case, $X_t$ and $Z_t$ share
the same limiting variance and the same limiting distribution. For
example, the distribution of $Z_t$ may be that of the Horvitz--Thompson
estimator (see \cite{hortho52}) under multistage sampling with simple
random without replacement sampling (SI) of PSUs, and the distribution
of $X_t$ may be that of the Hansen--Hurwitz estimator (see \cite
{hanhur43}) under multistage sampling and simple random with
replacement sampling (SIR) of PSUs.

In this paper, we derive asymptotic normality results for without-repla\-cement multistage designs, and we prove the consistency of a
with-replace\-ment bootstrap of PSUs for SI sampling at the first stage
when the sampling fraction tends to zero. Our framework and our
assumptions are defined in Section~\ref{secframe}. In Section~\ref{asympnor}, we first give an overview of asymptotic normality results
in survey sampling. We then state a central limit theorem for the
Horvitz--Thompson estimator in case of multistage sampling with
Bernoulli sampling (BE) of PSUs. The theorem follows from standard
assumptions and from the independent selections of PSUs. We generalize
to the multistage context a coupling algorithm by \cite{haj60} for the
joint selection of a BE sample and an SI sample. This is the main tool
to extend the central limit theorem to multistage sampling with SI
sampling of PSUs. We also prove the weak consistency of variance
estimators (see \cite{shatu95}, page~20), which enables to compute
normality-based confidence intervals with appropriate coverage. In
Section~\ref{grossecwrbootmult}, we consider the bootstrap for
multistage sampling. We introduce a new coupling algorithm between SI
sampling of PSUs and SIR sampling of PSUs. This is the main tool to
prove a long-standing issue; namely, that the so-called
with-replacement bootstrap of PSUs (see \cite{raowu88}) is consistent
in case of SI sampling of PSUs when the first-stage sampling fraction
tends to zero. This entails that Studentized bootstrap confidence
intervals are valid in such case, and that the bootstrap variance
estimators are consistent for smooth functions of totals. The
properties of a simplified variance estimator and of the bootstrap
procedure are evaluated in Section~\ref{secsimus} through a simulation
study. An application of the studied bootstrap method on the panel for
urban policy survey is presented in Section~\ref{applippv}. The proofs
of theorems are given in Section~\ref{secproofs}. Additional proofs
are given in the supplement \cite{cha15}.

\section{Framework} \label{secframe}

We consider a finite population $U$ consisting of $N$ sampling units
that may be represented by their labels, so that we may simply write
$U=\{1,\ldots,N\}$. The units are grouped inside $N_I$ nonoverlapping
sub-populations $u_1,\ldots,u_{N_I}$ called primary sampling units
(PSUs). We are interested in estimating the population total
\[
Y=\sum_{k \in U} y_k=\sum
_{u_i \in U_{I}} Y_i
\]
for some variable of interest $y$, where $Y_i=\sum_{k \in u_i} y_k$ is
the sub-total of the variable $y$ on the PSU $u_i$. We note $E(\cdot)$
and $V(\cdot)$ for the expectation and the variance of some estimator.
Also, we note $E_{\{X\}}(\cdot)$ and $V_{\{X\}}(\cdot)$ for the
expectation and variance conditionally on some random variable $X$.
Throughout the paper, we denote by $\hat{Y}_i$ an unbiased estimator of
$Y_i$, and by $V_i = V(\hat{Y}_i)$ its variance. Also, we denote by
$\hat{V}_i$ an unbiased estimator of $V_i$. In order to study the
asymptotic properties of the sampling designs and estimators that we
treat below, we consider the asymptotic framework of \cite{isaful82}.
We assume that the population $U$ belongs to a nested sequence $\{U_t\}
$ of finite populations with increasing sizes $N_t$, and that the
population vector of values $y_{Ut}=(y_{1t},\ldots,y_{Nt})^{\top}$
belongs to a sequence $\{y_{Ut}\}$ of $N_t$-vectors. For simplicity,
the index $t$ will be suppressed in what follows but all limiting
processes will be taken as $t \to\infty$.

In the population $U_I=\{u_1,\ldots,u_{N_I}\}$ of PSUs, a first-stage
sample $S_I$ is selected according to some sampling design $p_I(\cdot
)$. For clarity of exposition, we consider nonstratified sampling
designs for $p_I(\cdot)$, but the results may be easily extended to the
case of stratified first-stage sampling designs with a finite number of
strata, as is illustrated in Section~\ref{applippv}. If the PSU $u_i$
is selected in $S_I$, a second-stage sample $S_{i}$ is selected in
$u_i$ by means of some sampling design $p_i(\cdot|S_I)$. We assume
invariance of the second-stage designs: that is, the second stage of
sampling is independent of $S_I$ and we may simply write $p_i(\cdot
|S_I) = p_i(\cdot)$. Also, we assume that the second-stage designs are
independent from one PSU to another, conditionally on $S_I$. This
implies that
%
\begin{eqnarray}
\label{invindassumption} && \operatorname{Pr}\biggl(\bigcup_{u_i \in S_I}
\{S_i=s_{i}\} \Big| S_I\biggr) = \prod
_{u_i \in S_I} p_i(s_{i}|S_I) =
\prod_{u_i \in S_I} p_i(s_{i})
\end{eqnarray}
for any set of samples $s_{i} \subset u_i, i=1,\ldots,N_I$ (see \cite
{sarswewre92}, Chapter~4). The second-stage sampling designs
$p_i(\cdot)$ are left arbitrary. For example, they may involve censuses
inside some PSUs (which means cluster sampling), or additional stages
of sampling.

We will make use of the following assumptions:
\begin{longlist}[H3:]
\item[H1:] $N_I \displaystyle{\mathop{\longrightarrow}_{t\to\infty}} \infty$ and $n_I
\displaystyle{\mathop{\longrightarrow}_{t\to\infty}} \infty$. Also, $f_I=n_I/N_I
\displaystyle{\mathop{\longrightarrow}_{t\to\infty}} f \in[0,1[\,.$
\item[H2:] There exists $\delta>0$ and some constant $C_1$ such that
\begin{eqnarray*}
&& N_I^{-1} \sum_{u_i \in U_I} E|
\hat{Y}_i|^{2+\delta} < C_1.
\end{eqnarray*}
\item[H3:] There exists some constant $C_2$ such that $N_I^{-1} \sum_{u_i \in U_I} E(\hat{V}_i^{2}) < C_2$.
\item[H4:] There exists some constant $C_3>0$ such that
\begin{eqnarray*}
&& N_I^{-1} \sum_{u_i \in U_I}
(Y_i-\mu_Y)^2 > C_3 \qquad
\mbox{where } \mu_Y=N_I^{-1} Y.
\end{eqnarray*}
\end{longlist}

It is assumed in (H1) that a large number $n_I$ of PSUs is selected.
The assumption (H2) implies that the sequence of $\{Y_i\}_{u_i \in
U_I}$ has bounded moments of order $2+\delta$ and that the sequence of
$\{V(\hat{Y}_i)\}_{u_i \in U_I}$ has a bounded first moment. This
assumption requires in particular that the numbers of SSUs within PSUs
remain bounded. When we establish the mean square consistency of
variance estimators, assumption (H2) is strengthened by considering
$\delta=2$, which implies that the sequence of $\{Y_i\}_{u_i \in U_I}$
has bounded moments of order $4$. Assumptions (H2) and (H3) are
sufficient to have a weakly consistent variance estimator for further
stages of sampling. In this regard, assumption (H3) can be relaxed\vspace*{1pt} when
$f_I \displaystyle{\mathop{\longrightarrow}_{t\to\infty}} 0$ (see Section~\ref{secvestsrs}). Assumption (H4) requires that the dispersion between
PSUs does not vanish. This is a sufficient condition for the
first-stage sampling variance of the Horvitz--Thompson estimator to
have the usual order $O(N_I^{2} n_I^{-1})$, for the sampling designs
that we consider in this article.

\section{Asymptotic normality for multistage sampling} \label{asympnor}

Unbiased estimators for population totals such as the Horvitz--Thompson
estimator are well known; see \cite{hortho52} and \cite{nar51}.
Several results of asymptotic normality have been proved for specific
one-stage sampling designs; see, for example, \cite{haj60,haj61} for
simple random sampling without replacement, \cite{haj64} for rejective
sampling, \cite{ros72,sen80} and \cite{gor83} for successive
sampling, and \cite{ohl86} for the Rao--Hartley--Cochran procedure
proposed by \cite{raoharcoc62}. Br\"{a}nd\'en and Jonasson \cite{brajon12} state a
central limit theorem for the class of sampling algorithms satisfying
the strongly Rayleigh property, which includes Sampford sampling,
Pareto sampling and ordered pivotal sampling (see \cite{cha12}). Chen
and Rao \cite{cherao97} prove asymptotic normality for a class of estimators
under two-phase sampling designs; see also \cite{saewel13}. However,
asymptotic normality of estimators resulting from multistage samples
has not been much considered in the \mbox{literature}; two notable exceptions
are \cite{krerao81} for stratified multistage designs and
with-replacement sampling at the first-stage, and \cite{ohl89} who
states a martingale central limit theorem for a general two-stage
sampling design.

In this section, we confine our attention to Horvitz--Thompson
estimators for multistage sampling with BE sampling or SI sampling of
PSUs. The central limit Theorems \ref{theotclber} and \ref
{theotclsrs} are easily extended to cover smooth functions of totals
by using the delta method (see \cite{shatu95}, Appendix~A2).

\subsection{Bernoulli sampling of PSUs} \label{secber}

We first consider the case when a first-stage sample $S_I^B$ is
selected in $U_I$ by means of Bernoulli sampling (BE) with expected
size $n_I$, which we note as $S_I^B \sim \operatorname{BE}(U_I;n_I)$. The PSUs are
independently selected in $S_I^B$ with inclusion probabilities
$f_I=n_I/N_I$, and the size $n_{I}^B$ of $S_I^B$ is random. The
Horvitz--Thompson estimator
%
\begin{eqnarray}
\label{esthtber} \hat{Y}_B & = & \frac{N_I}{n_I} \sum
_{u_i \in U_{I}} I_i^B \hat{Y}_{i} =
\frac{N_I}{n_I} \sum_{u_i \in S_{I}^B} \hat{Y}_{i}
\end{eqnarray}
is unbiased for $Y$, with $I_i^B$ the sample membership indicator for
the PSU $u_i$ in the sample $S_I^B$. The variance of $\hat{Y}_B$ is
%
\begin{eqnarray}
\label{varhtber} V(\hat{Y}_B) &=& \frac{N_I^2}{n_I} \biggl
\{(1-f_I) \frac{1}{N_I} \sum_{u_i \in U_I}
Y_i^2 + \frac{1}{N_I} \sum
_{u_i \in U_I} V_i \biggr\},
\end{eqnarray}
where $V_i=V(\hat{Y}_i)$. We consider the variance estimator
%
\begin{eqnarray}
\label{estvarhtber} v_B(\hat{Y}_B) & = & \frac{N_I^2}{n_I}
\biggl(\frac{1-f_{I}}{n_I^B} \sum_{u_i \in S_I^B}
\hat{Y}_i^2 + \frac{f_I}{n_I^B} \sum
_{u_i \in S_I^B} \hat{V}_i\biggr)
\end{eqnarray}
if $n_I^B>0$ and $v_B(\hat{Y}_B)=0$ if $n_I^B=0$, with $\hat{V}_i$ an
unbiased estimator of $V_i$. Conditionally on $n_I^B$, $S_I^B$ may be
seen as an SI sample of size $n_I^B$ selected in $U_I$. It follows that
%
\begin{eqnarray}
E_{\{n_I^B\}} \biggl(\frac{1}{n_I^B} \sum_{u_i \in S_I^B}
\hat {Y}_i^2 \biggr) & = & \frac{1}{N_I} \sum
_{u_i \in U_I} \bigl(Y_i^2+V_i
\bigr),
\nonumber
\\[-8pt]
\label{estnIbber}
\\[-8pt]
\nonumber
E_{\{n_I^B\}} \biggl(\frac{1}{n_I^B} \sum
_{u_i \in S_I^B} \hat{V}_i \biggr) & = & \frac{1}{N_I}
\sum_{u_i \in U_I} V_i,
\end{eqnarray}
and $v_B(\hat{Y}_B)$ is unbiased for $V(\hat{Y}_B)$ conditionally on
$n_I^B=k>0$.

\begin{thm} \label{theotclber}
Assume that \textup{(H1)} and \textup{(H2)} hold. Then the Horvitz--Thompson estimator
$\hat{Y}_B = N_I n_I^{-1} \sum_{u_i \in S_{I}^B} \hat{Y}_{i}$ is
asymptotically normally distributed, that is,
%
\begin{eqnarray}
\label{tclber} && \bigl\{V(\hat{Y}_B)\bigr\}^{-0.5} (
\hat{Y}_B-Y) \displaystyle{\mathop{\longrightarrow}_{\mathcal{L}}}
\mathcal{N}(0,1),
\end{eqnarray}
where $\displaystyle{\mathop{\longrightarrow}_{\mathcal{L}}}$ stands\vspace*{-2.5pt}
for the convergence in distribution. Assume further that \textup{(H2)} holds
with $\delta=2$ and that \textup{(H3)} holds.\vspace*{1pt} Then $v_B(\hat{Y}_B)$ is
mean-square consistent for $V(\hat{Y}_B)$ conditionally on $n_I^B>0$,
that is,
%
\begin{eqnarray}
\label{estvarconsber0} E_{\{n_I^B>0\}} \bigl\llvert N_I^{-2}
n_I \bigl\{v_B(\hat{Y}_B)-V(\hat
{Y}_B) \bigr\}\bigr\rrvert ^2 & \displaystyle\mathop{
\longrightarrow}_{t\to\infty} & 0.
\end{eqnarray}
Also, $v_B(\hat{Y}_B)$ is mean-square consistent unconditionally:
%
\begin{eqnarray}
\label{estvarconsber} && E \bigl\llvert N_I^{-2} n_I
\bigl\{v_B(\hat{Y}_B)-V(\hat{Y}_B) \bigr\}
\bigr\rrvert ^2 \displaystyle{\mathop{\longrightarrow}_{t\to\infty}}
0.
\end{eqnarray}
\end{thm}

If $n_I^B>0$, we define $T_B \equiv\{v_B(\hat{Y}_B)\}^{-0.5} (\hat
{Y}_B-Y)$. It follows by the mean-square consistency of $v_B(\hat
{Y}_B)$ in (\ref{estvarconsber0}) that under assumption (H4),
$v_B(\hat{Y}_B)$ is weakly consistent for $V(\hat{Y}_B)$, namely
%
\begin{eqnarray}
\label{weakconsber} \bigl\{V(\hat{Y}_B)\bigr\}^{-1}
v_B(\hat{Y}_B) & \displaystyle\mathop{
\longrightarrow}_{{\mathrm{Pr}_{\{n_I^B>0\}}}} & 1,
\end{eqnarray}
where $\displaystyle\mathop{\longrightarrow}_{\mathrm{Pr}_{\{n_I^B>0\}}}$
stands for the convergence in
probability, conditionally on \mbox{$n_I^B>0$}. It follows by
(\ref{weakconsber}) and by the central limit theorem in (\ref{tclber})
that the pivotal quantity $T_B$ has a limiting standard normal
distribution. An approximate\vspace*{1.5pt} two-sided $100(1-2\alpha) \% $
confidence
interval for $Y$ is thus given by $ [\hat{Y}_B \pm u_{1-\alpha} \{
v_B(\hat{Y}_B)\}^{0.5} ]$, with $u_{1-\alpha}$ the quantile of
order $1-\alpha$ of the standard normal distribution.

\subsection{Without replacement simple random sampling of PSUs} \label{secsrs}

We consider the case when a first-stage sample $S_I$ is selected in
$U_I$ by means of simple random sampling without replacement (SI) of
size $n_I$, which we note as $S_I \sim \operatorname{SI}(U_I;n_I)$. The
Horvitz--Thompson estimator is
%
\begin{eqnarray}
\label{esthtsrs} \hat{Y} & = & \frac{N_I}{n_I} \sum
_{u_i \in U_{I}} I_i \hat{Y}_{i} =
\frac{N_I}{n_I} \sum_{u_i \in S_{I}} \hat{Y}_{i},
\end{eqnarray}
with $I_i$ the sample membership indicator for the PSU $u_i$ in the
sample $S_I$. We may alternatively rewrite the Horvitz--Thompson
estimator as
%
\begin{eqnarray}
\label{esthtsrs2} && \hat{Y} = N_I \bar{Z} \qquad \mbox{with }
\bar{Z}=\frac{1}{n_I} \sum_{j=1}^{n_I}
Z_j,
\end{eqnarray}
where the sample $S_I$ of PSUs is obtained by drawing $n_I$ times
without replacement one PSU in $U_I$, and where $Z_j$ stands for the
estimator of the total for the PSU selected at the $j$th draw. The
variance of $\hat{Y}$ is
%
\begin{eqnarray}
\label{varhtsrss}
&& V(\hat{Y}) = \frac{N_I^2}{n_I} \biggl\{(1-f_{I})
S_{Y,U_{I}}^2 + \frac
{1}{N_I} \sum
_{u_i \in U_I} V_i \biggr\},
\end{eqnarray}
with $S_{Y,U_{I}}^2=(N_{I}-1)^{-1} \sum_{u_i \in U_{I}} (Y_i-\mu
_{Y})^2$ the population dispersion of the sub-totals $Y_i$. Under (H1)
and (H2), $\hat{Y}$ is mean-square consistent for $Y$ in the sense that
%
\begin{eqnarray}
\label{consZbar} E\bigl\{N_I^{-1}(\hat{Y}-Y)\bigr
\}^2 & \displaystyle\mathop{\longrightarrow}_{t\to\infty} & 0.
\end{eqnarray}
This\vspace*{-1.5pt} implies that $N_I^{-1} (\hat{Y}-Y) \displaystyle\mathop{\longrightarrow}_{\mathrm{Pr}} 0$ where
$\displaystyle\mathop{\longrightarrow}_{\mathrm{Pr}}$ stands for the convergence in probability.

Hajek (1960) proposed a coupling procedure to draw simultaneously a BE
sample and an SI sample. This procedure is adapted in Algorithm \ref
{algo1} to the context of multistage sampling, and Proposition~\ref
{prop2} below generalizes the Lemma~2.1 in~\cite{haj60}.

\renewcommand{\thealgorithm}{\arabic{section}.\arabic{algorithm}}
\begin{algorithm}[t]
\caption{A coupling procedure for Bernoulli sampling of PSUs and simple
random sampling without replacement of PSUs} \label{algo1}
\begin{enumerate}
\item Draw the sample $S_I^B \sim \operatorname{BE}(U_I;n_I)$. Denote by $n_I^B$ the
(random) size of $S_I^B$.
\item Draw the sample $S_I$ as follows:
\begin{itemize}
\item if\vspace*{1pt} $n_I^B=n_I$, take $S_I=S_I^B$;
\item if\vspace*{1pt} $n_I^B<n_I$, draw $S_I^{+} \sim \operatorname{SI}(U_I \setminus
S_I^B;n_I-n_I^B)$ and take $S_I=S_I^B \cup S_I^{+}$;
\item if\vspace*{1pt} $n_I^B>n_I$, draw $S_I^{+} \sim \operatorname{SI}(S_I^B;n_I^B-n_I)$ and take
$S_I=S_I^B \setminus S_I^{+}$.
\end{itemize}
\item For any PSU $u_i$:
\begin{itemize}
\item if $u_i \in S_{I}^B \cap S_I$, select the same second-stage
sample $S_i$ for both $\hat{Y}_B$ and $\hat{Y}$;
\item if $u_i \in S_{I}^B \setminus S_I$, select a second-stage sample
$S_i$ for $\hat{Y}_B$;
\item if $u_i \in S_{I} \setminus S_I^B$, select a second-stage sample
$S_i$ for $\hat{Y}$.
\end{itemize}
\end{enumerate}
\end{algorithm}

\begin{prop} \label{prop2}
Assume that the samples $S_I^B$ and $S_I$ are selected according to
Algorithm \ref{algo1}. We note $\Delta_2 \equiv\sum_{u_i \in S_I}
(\hat{Y}_i-\mu_Y) - \sum_{u_i \in S_I^B} (\hat{Y}_i-\mu_Y)$. Then
%
\begin{eqnarray}
\label{inegsrsber} && \frac{E \{\Delta_2 \}^2}{V \{\sum_{u_i \in S_I^B} (\hat
{Y}_i-\mu_Y) \}} \leq \biggl\{\frac{1}{n_I}+
\frac{1}{N_I-n_I} \biggr\}^{0.5}.
\end{eqnarray}
\end{prop}

The result in Proposition~\ref{prop2} can be easily generalized to the
multivariate case: if $y_k=(y_{1k},\ldots,y_{qk})^{\top}$ denotes the
value taken for unit $k$ by some $q$-vector of interest $y$, we have
\begin{eqnarray*}
&& V \{\Delta_2 \} \leq \biggl\{\frac{1}{n_I}+
\frac
{1}{N_I-n_I} \biggr\}^{0.5} V \biggl\{\sum
_{u_i \in S_I^B} (\hat{Y}_i-\mu _Y) \biggr\},
\end{eqnarray*}
where for symmetric matrices $A$ and $B$ of size $q$, $A \leq B$ means
that $B-A$ is nonnegative definite.

\begin{thm} \label{theotclsrs}
Assume\vspace*{1pt} that \textup{(H1)} and \textup{(H2)} hold. Then the Horvitz--Thompson estimator
$\displaystyle\hat{Y} = N_I n_I^{-1} \sum_{u_i \in S_{I}} \hat{Y}_{i}$
is asymptotically normally distributed, that is,
%
\begin{eqnarray}
\label{tclsrs} \bigl\{V(\hat{Y})\bigr\}^{-0.5}(\hat{Y}-Y) &
\displaystyle\mathop{\longrightarrow}_{\mathcal{L}} & \mathcal{N}(0,1).
\end{eqnarray}
\end{thm}

\subsection{Variance estimation for SI sampling of PSUs} \label{secvestsrs}

We first consider the usual, unbiased variance estimator for $\hat{Y}$:
%
\begin{eqnarray}
&& v(\hat{Y}) = \frac{N_I^2}{n_I} \biggl\{(1-f_{I})
s_Z^2 + \frac{1}{N_I} \sum
_{u_i \in S_I} \hat{V}_i \biggr\}
\nonumber
\\[-8pt]
\label{estvarhtsrs}
\\[-8pt]
\eqntext{\displaystyle \mbox{with } s_Z^2=
\frac{1}{n_I-1} \sum_{j=1}^{n_I}
(Z_j-\bar{Z})^2.}
\end{eqnarray}

\begin{prop} \label{consestvarsrs}
Assume that \textup{(H1)} and \textup{(H3)} hold, and that \textup{(H2)} holds with $\delta=2$.
Then $v(\hat{Y})$ is mean-square consistent for $V(\hat{Y})$:
%
\begin{eqnarray}
\label{estvarconssrs} E \bigl\llvert N_I^{-2} n_I
\bigl\{v(\hat{Y})-V(\hat{Y}) \bigr\}\bigr\rrvert ^2 & \displaystyle
\mathop{\longrightarrow}_{t\to\infty} & 0.
\end{eqnarray}
\end{prop}

It follows by Proposition~\ref{consestvarsrs} that under the
assumption (H4), $v(\hat{Y})$ is weakly consistent for $V(\hat{Y})$, namely
%
\begin{eqnarray}
\label{weakconssrs} && \bigl\{V(\hat{Y})\bigr\}^{-1} v(\hat{Y})
\displaystyle\mathop{\longrightarrow}_{\mathrm{Pr}} 1.
\end{eqnarray}
From the central limit theorem in (\ref{tclsrs}), $T \equiv\{v(\hat
{Y})\}^{-0.5} (\hat{Y}-Y)$ has a limiting standard normal distribution.
Therefore, an approximate two-sided $100(1-2\alpha) \% $ confidence
interval for $Y$ is given by
%
\begin{eqnarray}
\label{icsrs} && \bigl[\hat{Y} \pm u_{1-\alpha} \bigl\{v(\hat{Y})\bigr
\}^{0.5} \bigr].
\end{eqnarray}

In proving Proposition~\ref{consestvarsrs}, assumption (H3) is
needed, requiring that an unbiased variance estimator $\hat{V}_i$ can
be computed inside PSUs. This assumption may be cumbersome,
particularly if the sampling design implies additional stages of
sampling inside PSUs. It is thus desirable to provide simplified
variance estimators which do not require assumption (H3) while
remaining consistent. We are able to do so in the particular important
case when the first-stage sampling rate tends to zero. A simplified
variance estimator (see \cite{sarswewre92}, equation (4.6.1)) is
obtained by simply dropping the term involving the variance estimators
inside PSUs $\hat{V}_i$. This leads to
%
\begin{eqnarray}
\label{estvarhtsrssimp1} && v_{\mathrm{SIMP}}(\hat{Y}) = \frac{N_I^2}{n_I}
(1-f_{I}) s_Z^2.
\end{eqnarray}

\begin{prop} \label{propvarsimp1}
Assume that \textup{(H1)} holds, and that \textup{(H2)} holds with $\delta=2$. Assume
that $f_I  \displaystyle\mathop{\longrightarrow}_{t\to\infty} 0$. Then $v_{\mathrm{SIMP}}(\hat
{Y})$ is mean-square consistent for $V(\hat{Y})$:
%
\begin{eqnarray}
\label{consvarsimp1} E \bigl\llvert N_I^{-2} n_I
\bigl\{v_{\mathrm{SIMP}}(\hat{Y})-V(\hat{Y}) \bigr\}\bigr\rrvert ^2 &
\displaystyle\mathop{\longrightarrow}_{t\to\infty} & 0.
\end{eqnarray}
\end{prop}

The proof (omitted) follows from the fact that when $f_I \displaystyle\mathop{\longrightarrow}_{t\to \infty} 0$, $V(\hat{Y})$ is asymptotically equivalent to
%
\begin{eqnarray}
\label{varhtsrs} && V_{\mathrm{app}}(\hat{Y}) = \frac{N_I^2}{n_I}
(1-f_{I}) \biggl\{S_{Y,U_{I}}^2 + \frac{1}{N_I}
\sum_{u_i \in U_I} V_i \biggr\}
\end{eqnarray}
under assumption (H2). It is easily seen from equation (\ref
{estvarhtsrs}) that $v_{\mathrm{SIMP}}(\hat{Y})$ tends to underestimate the
true variance, with a bias equal to $-\sum_{u_i \in U_I} V_i$. An
alternative simplified estimator is obtained by estimating the variance
as if the PSUs were selected with replacement [see equation (\ref{estvarwrhh})]. This leads to the second simplified variance estimator
%
\begin{eqnarray}
\label{estvarhtsrssimp2} && v_{\mathrm{WR}}(\hat{Y}) = \frac{N_I^2}{n_I}
s_Z^2.
\end{eqnarray}
It is easily shown that $v_{\mathrm{WR}}(\hat{Y})$ tends to overestimate the
true variance, with a bias equal to $N_I S_{Y,U_I}^2$. Under the
conditions of Proposition~\ref{propvarsimp1}, $v_{\mathrm{WR}}(\hat{Y})$ is
also mean-square consistent for the true variance since it only differs
from $v_{\mathrm{SIMP}}(\hat{Y})$ with the factor $(1-f_{I})$. Under the
additional assumption (H4), the variance estimators $v_{\mathrm{SIMP}}(\hat{Y})$
and $v_{\mathrm{WR}}(\hat{Y})$ are therefore weakly consistent for $V(\hat{Y})$.
When $f_I \displaystyle\mathop{\longrightarrow}_{t\to\infty} 0$, an approximate
two-sided $100(1-2\alpha) \% $ confidence interval for $Y$ is therefore
obtained from (\ref{icsrs}) by replacing $v(\hat{Y})$ with
$v_{\mathrm{SIMP}}(\hat{Y})$ or $v_{\mathrm{WR}}(\hat{Y})$.

\section{With-replacement bootstrap for multistage sampling} \label
{grossecwrbootmult}

The use of bootstrap techniques in survey sampling has been widely
studied in the literature. Most of them may be thought as particular
cases of the weighted bootstrap \cite{bercom97,anttil11,beapat12}; see also \cite{shatu95,davhin97,lah03} and \cite{davsar07} for detailed reviews.

Bootstrap for multistage sampling under without-replacement sampling of
PSUs has been considered, for example, in \cite{raowu88,raowuyue92,nigrao96,funsaisittoi06,pre09,linlurussit13},
among others. In this section, we consider the so-called
with-replacement bootstrap of PSUs (see \cite{raowu88}). This method
is suitable for with-replacement sampling of PSUs, and basic results
for such sampling designs are therefore reminded in Section~\ref{secwrsamp}. A~new coupling algorithm between SI sampling of PSUs and
SIR sampling of PSUs is given in Section~\ref{seccoupSIRSI}. This is
the main tool to study the bootstrap of PSUs for multistage sampling
with SI sampling of PSUs when the first-stage sampling fraction tends
to zero. In Section~\ref{secwrboot}, we prove that Studentized
bootstrap confidence intervals are valid. In Section~\ref{secvarestwrboot}, we prove that the bootstrap variance estimator is
consistent for smooth functions of means whenever it is consistent in
case of SIR sampling of PSUs.

\subsection{With replacement sampling of PSUs} \label{secwrsamp}

We consider the case when a first-stage sample $S_I^{\mathrm{WR}}$ is selected
in $U_I$ according to simple random sample with replacement (SIR) of
size $n_{I}$ inside $U_{I}$, which we note as $S_I^{\mathrm{WR}} \sim \operatorname{SIR}(U_I;n_I)$. Denote by $W_i$ the number of selections of the PSU
$u_i$ in $S_I^{\mathrm{WR}}$, and by $S_I^d$ of size $n_I^d$ the set of distinct
PSUs associated to $S_I^{\mathrm{WR}}$. Each time $j=1,\ldots,W_i$ that unit
$u_i$ is drawn in $S_I^{\mathrm{WR}}$, a second-stage sample $S_{i[j]}$ is
selected in $u_i$. The total $Y$ is unbiasedly estimated by the
Hansen--Hurwitz estimator
%
\begin{eqnarray}
\label{estwrhh} \hat{Y}_{\mathrm{WR}} & = & \sum_{u_i \in S_{I}^d}
\frac{1}{E(W_i)} \sum_{j=1}^{W_i}
\hat{Y}_{i[j]} = \frac{N_{I}}{n_{I}} \sum_{u_i \in S_I^d}
\sum_{j=1}^{W_i} \hat{Y}_{i[j]},
\end{eqnarray}
where $\hat{Y}_{i[j]}$ stands for an unbiased estimator of $Y_i$
computed on $S_{i[j]}$. We may alternatively rewrite the
Hansen--Hurwitz estimator as
%
\begin{eqnarray}
\label{estwrhh2} &&\hat{Y}_{\mathrm{WR}} = N_I \bar{X} \qquad
\mbox{with } \bar{X}=\frac{1}{n_I} \sum_{j=1}^{n_I}
X_j,
\end{eqnarray}
where the sample $S_I^{\mathrm{WR}}$ of PSUs is obtained by drawing $n_I$ times
with replacement one PSU in $U_I$ and where $X_j$ stands for the
estimator of the total for the PSU selected at the $j$th draw.

The variance of $\hat{Y}_{\mathrm{WR}}$ is
%
\begin{eqnarray}
\label{varwrhh} V (\hat{Y}_{\mathrm{WR}} ) & = & \frac{N_{I}^2}{n_{I}} \biggl\{
\frac
{N_{I}-1}{N_{I}} S_{Y,U_{I}}^2 + \frac{1}{N_{I}} \sum
_{u_i \in U_{I}} V_i \biggr\}.
\end{eqnarray}
Under (H1) and (H2), $\hat{Y}_{\mathrm{WR}}$ is mean-square consistent for $Y$
in the sense that
%
\begin{eqnarray}
\label{consXbar} E\bigl\{N_I^{-1}(\hat{Y}_{\mathrm{WR}}-Y)
\bigr\}^2 & \displaystyle\mathop{\longrightarrow}_{t\to\infty} & 0.
\end{eqnarray}
This implies that $N_I^{-1} (\hat{Y}_{\mathrm{WR}}-Y) \displaystyle \mathop{\longrightarrow}_{\mathrm{Pr}} 0$.

An unbiased variance estimator for $V (\hat{Y}_{\mathrm{WR}} )$ is
%
\begin{eqnarray}
\label{estvarwrhh} && v_{\mathrm{WR}} (\hat{Y}_{\mathrm{WR}} ) =
\frac{N_{I}^2}{n_{I}} s_X^2 \qquad \mbox{with }
s_X^2=\frac{1}{n_I-1} \sum
_{j=1}^{n_{I}} (X_j-\bar{X})^2.
\end{eqnarray}
The simple form of the variance estimator in (\ref{estvarwrhh}) is
primarily due to (\ref{estwrhh2}), where $\hat{Y}_{\mathrm{WR}}$ is written as
a sum of independent and identically distributed random variables (see
also \cite{sarswewre92}, page~151).

\begin{thm} \label{theotclwr}
Assume that \textup{(H1)} and \textup{(H2)} hold. Then the Hansen--Hurwitz estimator
$\hat{Y}_{\mathrm{WR}}=N_{I} n_{I}^{-1} \sum_{u_i \in S_I^d} \sum_{j=1}^{W_i}
\hat{Y}_{i[j]}$ is asymptotically normally distributed, that is,
%
\begin{eqnarray}
\label{tclwr} && \bigl\{V(\hat{Y}_{\mathrm{WR}})\bigr\}^{-0.5} (
\hat{Y}_{\mathrm{WR}}-Y) \displaystyle\mathop{\longrightarrow}_{\mathcal{L}}
\mathcal{N}(0,1).
\end{eqnarray}
Assume further that \textup{(H2)} holds with $\delta=2$. Then $v_{\mathrm{WR}} (\hat
{Y}_{\mathrm{WR}} )$ is mean-square consistent for $V (\hat
{Y}_{\mathrm{WR}} )$:
%
\begin{eqnarray}
\label{estvarconswr} E \bigl\llvert N_I^{-2} n_I
\bigl\{v_{\mathrm{WR}} (\hat{Y}_{\mathrm{WR}} )-V(\hat {Y}_{\mathrm{WR}}) \bigr
\}\bigr\rrvert ^2 & \displaystyle\mathop{\longrightarrow}_{t\to\infty}
& 0.
\end{eqnarray}
\end{thm}

In proving the consistency of $v_{\mathrm{WR}} (\hat{Y}_{\mathrm{WR}} )$,
assumption (H3) is not needed. In particular, unbiased variance
estimators $\hat{V}_i$ inside PSUs are not mandatory. This appealing
property leads to consider $v_{\mathrm{WR}}(\cdot)$ as a possible simplified
variance estimator when the PSUs are selected without replacement with
a first-stage sampling fraction tending to zero; see equation (\ref
{estvarhtsrssimp2}).

\subsection{A coupling procedure between SIR sampling of PSUs and SI
sampling of PSUs} \label{seccoupSIRSI}

The procedure is described in Algorithm \ref{algo2}. Conditionally on
$n_{I}^d$, the sample $S_{I}^d$ obtained in step~1 is by symmetry an SI
sample of size $n_{I}^d$ from $U_{I}$, which implies that $S_{I}^d \cup
S_{I}^c$ is an SI sample of size $n_{I}$ from $U_{I}$. Consequently,
this procedure leads to a sample $S_I$ drawn by means of SI sampling of PSUs.

\setcounter{algorithm}{0}
\begin{algorithm}
\caption{A coupling procedure for simple random sampling
with-replacement of PSUs and simple random sampling without replacement
of PSUs for multistage sampling} \label{algo2}
\begin{enumerate}
\item Draw the sample $S_I^{\mathrm{WR}} \sim \operatorname{SIR}(U_I;n_I)$. Denote by $S_{I}^d$
of (random) size $n_{I}^d$ the set of distinct PSUs in $S_I^{\mathrm{WR}}$.
\item Draw a complementary sample $S_I^{c} \sim \operatorname{SI}(U_I \setminus
S_I^d;n_I-n_I^d)$ and take $S_I=S_{I}^d \cup S_I^{c}$.
\item For any $u_i \in S_I^d$:
\begin{itemize}
\item Each time $j=1,\ldots,W_i$ that unit $u_i$ is drawn in
$S_I^{\mathrm{WR}}$, select a second-stage sample $S_{i[j]}$ with associated
estimator $\hat{Y}_{i[j]}$ for $\hat{Y}_{\mathrm{WR}}$.
\item Take $S_i=S_{i[1]}$ and $\hat{Y}_i=\hat{Y}_{i[1]}$ for $\hat{Y}$.
\end{itemize}
\item For any $u_i \in S_I^c$, select a second-stage sample $S_i$ with
associated estimator $\hat{Y}_i$ for $\hat{Y}$.
\end{enumerate}
\end{algorithm}

\begin{prop} \label{prop5}
Assume that the samples $S_I^{\mathrm{WR}}$ and $S_I$ are selected according to
Algorithm \ref{algo2}. Then
%
\begin{eqnarray}
\label{prop5eq1} \frac{E(\hat{Y}_{\mathrm{WR}}-\hat{Y})^2}{V(\hat{Y}_{\mathrm{WR}})} & \leq& \frac{n_I-1}{N_I-1}.
\end{eqnarray}
\end{prop}

The right bound in (\ref{prop5eq1}) is mainly of interest when $f_I
\displaystyle\mathop{\longrightarrow}_{t\to\infty} 0$. In this case, from the
trivial inequality $\frac{n_I-1}{N_I-1} \leq\frac{n_I}{N_I}$,
Algorithm \ref{algo2} may be used to select the samples $S_I^{\mathrm{WR}}$ and
$S_I$ so that the difference between $\hat{Y}_{\mathrm{WR}}$ and $\hat{Y}$ is
asymptotically negligible. A similar result holds for the dispersions
between the estimated totals inside PSUs, as stated in Proposition~\ref
{prop6} below.

\begin{prop} \label{prop6}
Assume that the samples $S_I^{\mathrm{WR}}$ and $S_I$ are selected according to
Algorithm \ref{algo2}. Assume that \textup{(H1)} and \textup{(H2)} hold, and that
$f_I \displaystyle\mathop{\longrightarrow}_{t\to\infty} 0$. Then
%
\begin{eqnarray}
\label{prop6eq1} E(\bar{Z}-\bar{X})^2 &=& o\bigl(n_I^{-1}
\bigr),
\\
\label{prop6eq2} E\bigl|s_Z^2-s_X^2\bigr|
& \displaystyle\mathop{\longrightarrow}_{t\to\infty}& 0,
\end{eqnarray}
where $\bar{X}$ and $s_X^2$ are defined in equations (\ref{estwrhh2})
and (\ref{estvarwrhh}), and $\bar{Z}$ and $s_Z^2$ are defined in
equations (\ref{esthtsrs2}) and (\ref{estvarhtsrs}),
\end{prop}

\subsection{With replacement bootstrap of PSUs} \label{secwrboot}

We consider the with-replace\-ment bootstrap of PSUs described in
\cite{raowu88}. Using the notation introduced in equation~(\ref
{esthtsrs2}), let $(Z_1,\ldots,Z_{n_I})^{\top}$ denote the sample of
estimators under SI sampling of PSUs. Also, let $(Z_1^*,\ldots
,Z_{m}^*)^{\top}$ be obtained by sampling $m$ times independently in
$(Z_1,\ldots,Z_{n_I})^{\top}$. Similarly, using the notation introduced
in equation~(\ref{estwrhh2}), let $(X_1,\ldots,X_{n_I})^{\top}$
denote the sample of estimators under SIR sampling of PSUs. Also, let
$(X_1^*,\ldots,X_{m}^*)^{\top}$ be obtained by sampling $m$ times
independently in $(X_1,\ldots,X_{n_I})^{\top}$.

We first demonstrate the bootstrap consistency. We note
\begin{eqnarray*}
\bar{Z}_m^* &=& \frac{1}{m} \sum
_{j=1}^m Z_j^* \quad \mbox{and}\quad
s_Z^{*2} = \frac{1}{m-1} \sum
_{j=1}^m \bigl(Z_j^*-
\bar{Z}_m^* \bigr)^2,
\\
\bar{X}_m^* &=& \frac{1}{m} \sum
_{j=1}^m X_j^* \quad \mbox{and} \quad
s_X^{*2} = \frac{1}{m-1} \sum
_{j=1}^m \bigl(X_j^*-
\bar{X}_m^* \bigr)^2.
\end{eqnarray*}

We proceed by showing that, using Algorithm \ref{algo2}, the samples
$S_I$ and $S_I^{\mathrm{WR}}$ can be drawn so that the pivotal statistics
%
\begin{eqnarray}
\label{pivotstat} && m^{0.5} \bigl(s_Z^{*}
\bigr)^{-1} \bigl(\bar{Z}_m^*-\bar{Z}\bigr) \quad \mbox{and}
\quad m^{0.5} \bigl(s_X^{*}\bigr)^{-1}
\bigl(\bar{X}_m^*-\bar{X}\bigr)
\end{eqnarray}
are close. More precisely, we make use of the Mallows metric (see \cite
{mal72} and \cite{bicfre81}), also known as the Wasserstein metric.
Let $1 \leq q < \infty$, and let $\alpha$ and $\beta$ denote two
distributions on $\mathbb{R}^s$ with finite moments of order $q$. Then
%
\begin{eqnarray}
\label{mallowmetric} d_q(\alpha,\beta) & = & \inf \bigl\{E\|X-Z
\|^q \bigr\}^{1/q},
\end{eqnarray}
where the infimum is taken over all couples $(X,Z)$ with marginal
distributions $\alpha$ and $\beta$. For two random vectors $X$ and $Z$,
we note $d_q(\alpha,\beta)$ for the $d_q$-distance between the
distributions of $X$ and $Z$. In what follows, we consider $q=1$ or $q=2$.

Let $D=(D_1,\ldots,D_{n_I})^{\top}$ be generated according to a
multinomial distribution with parameters $(m;n_I^{-1},\ldots
,n_I^{-1})$. The same multinomial weights $D$ are used in the selection
of $(Z_1^*,\ldots,Z_{m}^*)^{\top}$ and $(X_1^*,\ldots,X_{m}^*)^{\top}$,
so that we may write
%
\begin{eqnarray}
\label{estbootpond} && \bar{Z}_m^* = \frac{1}{m} \sum
_{j=1}^{n_I} D_j Z_j \quad
\mbox{and} \quad \bar{X}_m^* = \frac{1}{m} \sum
_{j=1}^{n_I} D_j X_j.
\end{eqnarray}

\begin{prop} \label{prop6b}
Assume that \textup{(H1)} and \textup{(H2)} hold. Assume that $f_I
\displaystyle\mathop{\longrightarrow}_{t\to\infty} 0$ and that
$m \displaystyle\mathop{\longrightarrow}_{t\to\infty} \infty$. Then
%
\begin{eqnarray}
\label{prop6beq1} && E \bigl(\bar{Z}_m^*-\bar{X}_m^*
\bigr)^2 = o\bigl(m^{-1}\bigr)+o\bigl(n_I^{-1}
\bigr).
\end{eqnarray}
\end{prop}

\begin{prop} \label{prop7}
Assume that \textup{(H1)} and \textup{(H2)} hold. Assume that $f_I
\displaystyle\mathop{\longrightarrow}_{t\to\infty } 0$ and that
$m \displaystyle\mathop{\longrightarrow}_{t\to\infty} \infty$. Then
%
\begin{eqnarray}
\label{d2ZmstarXmstar} d_2 \bigl[ m^{0.5} \bigl(
\bar{Z}_m^*-\bar{Z}\bigr),m^{0.5} \bigl(\bar{X}_m^*-
\bar {X}\bigr) \bigr] & \displaystyle\mathop{\longrightarrow}_{t\to\infty}& 0,
\\
\label{d2sZstarsXstar} d_1 \bigl[ s_Z^{*2},
s_X^{*2} \bigr] & \displaystyle\mathop{
\longrightarrow}_{t\to\infty}& 0,
\end{eqnarray}
where the distance $d_q(\cdot,\cdot)$ is defined in (\ref{mallowmetric}).
\end{prop}

From Proposition~\ref{prop7}, the pivotal statistics in (\ref
{pivotstat}) share the same limiting distribution. Theorem~\ref
{consboot} below follows from Theorem~2.1 of \cite{bicfre81}.

\begin{thm} \label{consboot}
Assume that \textup{(H1)} and \textup{(H2)} hold. Assume that
$f_I \displaystyle\mathop{\longrightarrow}_{t\to\infty } 0$ and that
$m \displaystyle\mathop{\longrightarrow}_{t\to\infty}\infty$. Then
%
\begin{eqnarray}
\label{bootpivquant} && m^{0.5} \bigl(s_Z^{*}
\bigr)^{-1} \bigl(\bar{Z}_m^*-\bar{Z}\bigr) \mathop{
\longrightarrow}_{\mathcal{L}} \mathcal{N}(0,1).
\end{eqnarray}
\end{thm}

Theorem~\ref{consboot} implies that the normality-based confidence
interval for $Y$ given in~(\ref{icsrs}) may be replaced by the
Studentized bootstrap confidence interval (see \cite{davhin97}, page~194)
%
\begin{eqnarray}
\label{icstubootsrs} && \bigl[\hat{Y} - u_{1-\alpha}^* \bigl\{v(\hat{Y})\bigr
\}^{0.5},\hat{Y} - u_{\alpha
}^* \bigl\{v(\hat{Y})\bigr
\}^{0.5} \bigr],
\end{eqnarray}
where the quantiles $u_{1-\alpha}$ and $u_{\alpha}$ of the normal
distribution are replaced by the corresponding quantiles $u_{1-\alpha
}^*$ and $u_{\alpha}^*$ of the bootstrap pivotal quantity in (\ref
{bootpivquant}). The simplified variance estimators $v_{\mathrm{SIMP}}(\hat
{Y})$ and $v_{\mathrm{WR}}(\hat{Y})$ can also be used in (\ref{icstubootsrs}).

\subsection{Bootstrap variance estimation for functions of totals}
\label{secvarestwrboot}

We now consider the case when $y_k=(y_{1k},\ldots,y_{qk})^{\top}$ is
multivariate, and denotes the value taken for unit $k$ by some
$q$-vector of interest $y$. We are interested in a parameter $\theta
=f(Y)$ for some function $f:\mathbb{R}^q \longrightarrow\mathbb{R}$.
Under SI sampling of PSUs, the plug-in estimator of $\theta$ is $\hat
{\theta} \equiv f(N_I \bar{Z})$. Under SIR sampling of PSUs, the
plug-in estimator of $\theta$ is $\hat{\theta}_{\mathrm{WR}} \equiv f(N_I \bar
{X})$. Also, we note $\hat{\theta}^* \equiv f(\bar{Z}_m^*)$ and $\hat
{\theta}_{\mathrm{WR}}^* \equiv f(\bar{X}_m^*)$ for the bootstrap estimators,
where $\bar{Z}_m^*$ and $\bar{X}_m^*$ are defined in (\ref
{estbootpond}). We consider the additional regularity assumptions:
\begin{itemize}[H6:]
\item[H5:] $f(\cdot)$ is homogeneous of degree $\beta\geq0$, in that
$f(ry)=r^{\beta}f(y)$ for any real $r>0$ and $q$-vector $y$. Also, $f$
is a differentiable function on $\mathbb{R}^q$ with bounded partial derivatives.
\item[H6:] There exists some constant $C_4>0$ such that $V(\hat{\theta
}_{\mathrm{WR}}) > C_4 N_I^{2\beta} n_I^{-1}$.
\end{itemize}

Assumption (H6) is similar to (H4), and requires the variance of the
plug-in estimator $\hat{\theta}_{\mathrm{WR}}$ to have the usual order
$O(N_I^{2\beta} n_I^{-1})$.

\begin{prop} \label{prop8}
Assume that the samples $S_I^{\mathrm{WR}}$ and $S_I$ are selected according to
Algorithm \ref{algo2}. Assume that assumptions \textup{(H1)}, \textup{(H2)} and \textup{(H5)}
hold. Assume that $f_I \displaystyle\mathop{\longrightarrow}_{t\to\infty} 0$. Then
%
\begin{eqnarray}
\label{prop8eq1} E\bigl(\|\bar{Z}-\bar{X}\|^2\bigr) &=& o
\bigl(n_I^{-1}\bigr),
\\
\label{prop8eq2} E(\hat{\theta}-\hat{\theta}_{\mathrm{WR}})^2 &=& o
\bigl(N_I^{2\beta} n_I^{-1}\bigr),
\end{eqnarray}
with $\|\cdot\|$ the Euclidean norm. Assume further that $m \displaystyle\mathop{\longrightarrow}_{t\to\infty} \infty$. Then
%
\begin{eqnarray}
\label{prop9eq1} E\bigl(\bigl\|\bar{Z}^*-\bar{X}^*\bigr\|^2\bigr) &=& o
\bigl(m^{-1}\bigr)+o\bigl(n_I^{-1}\bigr),
\\
\label{prop9eq2} E\bigl(\hat{\theta}^*-\hat{\theta}_{\mathrm{WR}}^*
\bigr)^2 &=& o\bigl(N_I^{2\beta} m^{-1}
\bigr)+o\bigl(N_I^{2\beta} n_I^{-1}
\bigr).
\end{eqnarray}
\end{prop}

\begin{prop} \label{prop9}
Assume that the samples $S_I^{\mathrm{WR}}$ and $S_I$ are selected according to
Algorithm \ref{algo2}. Assume that assumptions \textup{(H1)}, \textup{(H2)}, \textup{(H5)} and
\textup{(H6)} hold. Assume that $f_I \displaystyle\mathop{\longrightarrow}_{t\to\infty} 0$
and $m=O(n_I)$. Then
%
\begin{eqnarray}
&& \frac{V_{\{X\}}(\hat{\theta}_{\mathrm{WR}}^*)}{V(\hat{\theta}_{\mathrm{WR}})} \mathop{\longrightarrow}_{\mathrm{Pr}} 1 \quad
\mbox{implies} \quad \frac{V_{\{Z\}}(\hat{\theta
}^*)}{V(\hat{\theta})} \mathop{\longrightarrow}_{\mathrm{Pr}} 1,
\end{eqnarray}
with $V_{\{X\}}$ the variance conditionally on $X_1,\ldots,X_{n_I}$,
and similarly for $V_{\{Z\}}$.
\end{prop}

The proof or Proposition~\ref{prop8} follows from the regularity
assumptions on $f(\cdot)$ and from Propositions \ref{prop6} and \ref
{prop6b}. Proposition~\ref{prop9} implies that the with-replacement
bootstrap of PSUs provides consistent variance estimation for $\hat
{\theta}$ whenever it does so for $\hat{\theta}_{\mathrm{WR}}$. The regularity
assumption (H5) is somewhat strong, and may be weakened to
differentiability of $f(\cdot)$ on a compact set, under additional
assumptions on the vector of interest $y$ and on the second-stage
sampling weights.

\section{A simulation study} \label{secsimus}

We conducted a limited simulation study to investigate on the
performance of the variance estimators. We first generated 3 finite
populations, each with $N_I=2000$ PSUs. The number of SSUs inside PSUs
was generated so that the average number of SSUs per PSU was
approximately equal to $\bar{N}=40$, and so that the coefficient of
variation for the sizes $N_i$ of PSUs was equal to $0$ for population 1
(so that the PSUs are of equal size), approximately equal to $0.03$ for
population 2, and approximately equal to $0.06$ for population 3.

In each population, we generated for any PSU $u_i$ the value $\lambda_i
= \lambda+ \sigma~v_i$ with $\lambda=20$ and $\sigma=2$ for each
population, and the $v_i$'s were generated according to a normal
distribution with mean $0$ and variance $1$. For each SSU $k \in u_i$,
we generated three couples of values $(y_{1,k},y_{2,k})$,
$(y_{3,k},y_{4,k})$ and $(y_{5,k},y_{6,k})$ according to the model
%
\begin{eqnarray}
y_{2h-1,k} & = & \lambda_i + \bigl\{\rho_h^{-1}(1-
\rho_h)\bigr\}^{0.5} \sigma (\alpha \varepsilon_k+
\eta_k),
\\
y_{2h,k} & = & \lambda_i + \bigl\{\rho_h^{-1}(1-
\rho_h)\bigr\}^{0.5} \sigma (\alpha\varepsilon_k+
\nu_k),
\end{eqnarray}
for $h=1,\ldots,3$, where the values $\varepsilon_k$, $\eta_k$ and $\nu_k$
were generated according to a normal distribution with mean $0$ and
variance $1$. In each population, the parameter $\rho_h$ was chosen so
that the intra-cluster correlation coefficient was approximately equal
to $0.1$ for both variables $y_1$ and $y_2$, $0.2$ for both variables
$y_3$ and $y_4$, and $0.3$ for both variables $y_5$ and $y_6$. Also,
the parameter $\alpha$ was chosen so that the coefficient of
correlation between variables $y_{2h-1}$ and $y_{2h}$, $h=1,\ldots,3$,
was approximately equal to $0.60$.

From each population, we selected $B=1000$ two-stage samples. The
sample $S_I$ of PSUs was selected by means of SI sampling of size
$n_I=20,40,100$ or $200$. Inside each $u_i \in S_I$, the sample $S_i$
of SSUs was selected by means of systematic sampling of size $n_0=5$ or
$10$. Note that, due to the systematic sampling at the second stage,
the variance may not be unbiasedly estimated. Our objective is to
estimate the variance of the Horvitz--Thompson estimator of the totals
of the variables $y_1$, $y_3$ and $y_5$, by using the simplified
variance estimator $v_{\mathrm{SIMP}}(\hat{Y})$ in (\ref{estvarhtsrssimp1})
or the with-replacement bootstrap of PSUs. Also, our objective is to
estimate the variance of the substitution estimator for the ratios
%
\begin{eqnarray}
\label{ratio} R_h & = & (\mu_{y,2h})^{-1}
\mu_{y,2h-1}
\end{eqnarray}
with $\mu_{y,2h-1}=N^{-1} \sum_{k \in U} y_{2h-1,k}$ and $\mu
_{y,2h}=N^{-1} \sum_{k \in U} y_{2h,k}$, and the variance for the
substitution estimator for the coefficient of correlations
%
\begin{eqnarray}
\label{coefcorr} r_h & = & \frac{\sum_{k \in U} (y_{2h-1,k}-\mu_{y,2h-1})(y_{2h,k}-\mu
_{y,2h})}{\{\sum_{k \in U} (y_{2h-1,k}-\mu_{y,2h-1})^2 \sum_{k \in U}
(y_{2h,k}-\mu_{y,2h})^2\}^{0.5}},
\end{eqnarray}
for $h=1,\ldots,3$ by using the with-replacement bootstrap of PSUs. The
true variance was approximated from a separate simulation run of
$C=20{,}000$ samples.

As a measure of bias of a point estimator $\hat{\theta}$ of a parameter
$\theta$, we used the Monte Carlo percent relative bias (RB) given by
\[
\mathrm{RB}_{\mathrm{MC}}(\hat{\theta})=100 \frac{B^{-1} \sum_{b=1}^B \hat{\theta
}_{(b)} -\theta}{\theta},
\]
where $\hat{\theta}_{(b)}$ gives the value of the estimator for the
$b$th sample. As a measure of variance of an estimator $\hat{\theta,}$
we used the Monte Carlo percent relative stability (RS) given by
\[
\mathrm{RS}_{\mathrm{MC}}(\hat{\theta})=100 \frac{\{B^{-1} \sum_{b=1}^B (\hat{\theta
}_{(b)} -\theta)^2\}^{0.5}}{\theta}.
\]
When the simplified variance estimator in (\ref{estvarhtsrssimp1})
is used, we also assess the coverage of confidence intervals based on
asymptotic normality. When the with-replacement bootstrap of PSUs is
used, we assess the coverage of confidence intervals obtained by means
of the percentile method. We used a nominal one-tailed error rate of
2.5\% in each tail.
\begin{table}
\caption{Relative bias, relative stability and nominal one-tailed
error rates for the simplified variance estimator of the
Horvitz--Thompson estimator, and for the bootstrap for the estimation
of  a total for population~3} \label{tabsimres3a}
\begin{tabular*}{\tablewidth}{@{\extracolsep{\fill}}lcd{2.2}d{1.2}d{2.2}d{2.2}d{2.2}d{1.2}d{2.2}d{2.2}@{}}
\hline
& & \multicolumn{8}{c@{}}{\textbf{Simplified variance estimator} $\bolds{v_{\mathrm{SIMP}}(\hat
{Y})}$ \textbf{for} $\bolds{Y_{2h-1}}$} \\[-4pt]
&&\multicolumn{8}{l@{}}{\hrulefill} \\
& \multicolumn{1}{c}{$\bolds{n_0}$} & \multicolumn{4}{c}{$\mathbf{5}$} & \multicolumn{4}{c@{}}{$\textbf{10}$} \\
&& \multicolumn{4}{l}{\hrulefill}  & \multicolumn{4}{l@{}}{\hrulefill}\\
& \multicolumn{1}{c}{$\bolds{n_I}$} &
 \multicolumn{1}{c}{$\bolds{20}$} &
 \multicolumn{1}{c}{$\bolds{40}$} & \multicolumn{1}{c}{$\bolds{100}$} &
 \multicolumn{1}{c}{$\bolds{200}$} & \multicolumn{1}{c}{$\bolds{20}$} &
 \multicolumn{1}{c}{$\bolds{40}$} & \multicolumn{1}{c}{$\bolds{100}$} &
 \multicolumn{1}{c@{}}{$\bolds{200}$} \\
\hline
$\rho=0.1$ & RB & -0.02 & 0.01 & -0.03 & -0.04 & 0.00 & 0.01 & -0.01 &
-0.03 \\
& RS & 0.31 & 0.23 & 0.14 & 0.10 & 0.31 & 0.22 & 0.14 & 0.10 \\
& L & 2.9 & 3.0 & 2.0 & 2.3 & 2.4 & 2.5 & 2.4 & 2.3 \\
& U & 3.7 & 3.0 & 2.3 & 2.1 & 3.0 & 2.9 & 2.8 & 3.9 \\
& L${}+{}$U & 6.6 & 6.0 & 4.3 & 4.4 & 5.4 & 5.4 & 5.2 & 6.2 \\[3pt]
$\rho=0.2$ & RB & -0.04 & 0.00 & -0.01 & -0.04 & -0.01 & 0.01 & -0.01
& -0.01 \\
& RS & 0.31 & 0.21 & 0.14 & 0.10 & 0.32 & 0.23 & 0.14 & 0.09 \\
& L & 3.6 & 3.3 & 2.4 & 2.1 & 3.2 & 3.7 & 1.9 & 3.0 \\
& U & 3.6 & 3.3 & 2.0 & 1.9 & 2.7 & 3.1 & 2.1 & 2.6 \\
& L${}+{}$U & 7.2 & 6.6 & 4.4 & 4.0 & 5.9 & 6.8 & 4.0 & 5.6 \\[3pt]
$\rho=0.3$ & RB & -0.03 & 0.01 & -0.01 & -0.02 & 0.00 & 0.02 & 0.00 &
-0.02 \\
& RS & 0.31 & 0.22 & 0.13 & 0.09 & 0.33 & 0.22 & 0.14 & 0.09 \\
& L & 3.0 & 2.9 & 2.0 & 3.1 & 4.3 & 3.0 & 2.5 & 2.6 \\
& U & 3.1 & 2.6 & 2.1 & 1.9 & 3.1 & 3.5 & 2.2 & 3.7 \\
& L+U & 6.1 & 5.5 & 4.1 & 5.0 & 7.4 & 6.5 & 4.7 & 6.3 \\[6pt]
& & \multicolumn{8}{c@{}}{\textbf{Bootstrap of PSUs for} $\bolds{Y_{2h-1}}$} \\[-4pt]
&& \multicolumn{8}{l@{}}{\hrulefill}\\
$\rho=0.1$ & RB & -0.01 & 0.03 & 0.02 & 0.06 & 0.02 & 0.03 & 0.04 &
0.08 \\
& RS & 0.31 & 0.24 & 0.15 & 0.13 & 0.32 & 0.23 & 0.16 & 0.14 \\
& L & 3.3 & 2.7 & 2.4 & 1.7 & 2.7 & 2.4 & 2.2 & 1.8 \\
& U & 3.8 & 2.9 & 2.3 & 1.6 & 3.1 & 2.9 & 2.4 & 2.6 \\
& L${}+{}$U & 7.1 & 5.6 & 4.7 & 3.3 & 5.8 & 5.3 & 4.6 & 4.4 \\[3pt]
$\rho=0.2$ & RB & -0.03 & 0.02 & 0.04 & 0.07 & 0.01 & 0.03 & 0.04 &
0.10 \\
& RS & 0.32 & 0.22 & 0.16 & 0.13 & 0.33 & 0.24 & 0.16 & 0.15 \\
& L & 3.6 & 3.5 & 2.1 & 1.8 & 3.1 & 3.3 & 2.1 & 2.4 \\
& U & 3.5 & 3.3 & 1.9 & 1.4 & 2.8 & 3.0 & 1.8 & 1.8 \\
& L${}+{}$U & 7.1 & 6.8 & 4.0 & 3.2 & 5.9 & 6.3 & 3.9 & 4.2 \\[3pt]
$\rho=0.3$ & RB & -0.02 & 0.02 & 0.04 & 0.09 & 0.02 & 0.04 & 0.06 &
0.09 \\
& RS & 0.32 & 0.23 & 0.15 & 0.14 & 0.34 & 0.24 & 0.16 & 0.15 \\
& L & 2.9 & 3.1 & 1.8 & 2.2 & 4.3 & 2.9 & 2.4 & 2.3 \\
& U & 3.3 & 2.9 & 1.8 & 1.8 & 3.4 & 3.5 & 2.4 & 2.4 \\
& L${}+{}$U & 6.2 & 6.0 & 3.6 & 4.0 & 7.7 & 6.4 & 4.8 & 4.7 \\
\hline
\end{tabular*}
\end{table}
\begin{table}
\caption{Relative bias, relative stability and nominal one-tailed
error rates for the bootstrap for the estimation of a ratio and a
coefficient of correlation for population 3} \label{tabsimres3b}
\begin{tabular*}{\tablewidth}{@{\extracolsep{\fill}}lcd{2.2}d{1.2}d{2.2}d{2.2}d{2.2}d{1.2}d{1.2}d{1.2}@{}}
\hline
& & \multicolumn{8}{c@{}}{\textbf{Bootstrap of PSUs for $\bolds{R_h}$}} \\[-4pt]
&& \multicolumn{8}{l@{}}{\hrulefill}\\
& \multicolumn{1}{c}{$\bolds{n_0}$} & \multicolumn{4}{c}{$\bolds{5}$} & \multicolumn{4}{c@{}}{$\bolds{10}$}
\\[-4pt]
&& \multicolumn{4}{l}{\hrulefill} & \multicolumn{4}{l@{}}{\hrulefill}\\
& \multicolumn{1}{c}{$\bolds{n_I}$} &
\multicolumn{1}{c}{$\bolds{20}$} &
\multicolumn{1}{c}{$\bolds{40}$} &
\multicolumn{1}{c}{$\bolds{100}$} &
\multicolumn{1}{c}{$\bolds{200}$} &
\multicolumn{1}{c}{$\bolds{20}$} &
\multicolumn{1}{c}{$\bolds{40}$} &
\multicolumn{1}{c}{$\bolds{100}$} &
\multicolumn{1}{c@{}}{$\bolds{200}$} \\
\hline
$\rho=0.1$ & RB & 0.03 & 0.01 & 0.03 & 0.01 & 0.01 & 0.00 & 0.02 &
0.03 \\
& RS & 0.33 & 0.24 & 0.16 & 0.11 & 0.34 & 0.24 & 0.16 & 0.12 \\
& L & 2.6 & 3.5 & 2.5 & 2.6 & 3.8 & 2.9 & 3.0 & 4.0 \\
& U & 3.3 & 3.7 & 3.3 & 3.1 & 2.6 & 3.1 & 2.5 & 2.5 \\
& L${}+{}$U & 5.9 & 7.2 & 5.8 & 5.7 & 6.4 & 6.0 & 5.5 & 6.5 \\[3pt]
$\rho=0.1$ & RB & 0.01 & 0.01 & 0.00 & 0.02 & 0.01 & 0.01 & 0.01 &
0.03 \\
& RS & 0.34 & 0.23 & 0.15 & 0.12 & 0.34 & 0.23 & 0.14 & 0.12 \\
& L & 3.2 & 2.4 & 2.7 & 2.3 & 2.6 & 2.4 & 2.2 & 2.9 \\
& U & 3.0 & 2.4 & 2.2 & 2.2 & 4.4 & 4.0 & 1.8 & 2.5 \\
& L${}+{}$U & 6.2 & 4.8 & 4.9 & 4.5 & 7.0 & 6.4 & 4.0 & 5.4 \\[3pt]
$\rho=0.1$ & RB & -0.02 & 0.02 & 0.02 & 0.03 & 0.01 & 0.02 & 0.02 &
0.02 \\
& RS & 0.32 & 0.24 & 0.15 & 0.12 & 0.35 & 0.24 & 0.15 & 0.12 \\
& L & 3.5 & 2.8 & 2.7 & 3.6 & 2.2 & 2.4 & 2.9 & 3.2 \\
& U & 3.8 & 3.5 & 1.9 & 2.4 & 3.3 & 2.2 & 2.4 & 2.4 \\
& L${}+{}$U & 7.3 & 6.3 & 4.6 & 6.0 & 5.5 & 4.6 & 5.3 & 5.6 \\[6pt]
& & \multicolumn{8}{c@{}}{\textbf{Bootstrap of PSUs for} $\bolds{r_h}$} \\[-4pt]
&& \multicolumn{8}{l@{}}{\hrulefill}\\
$\rho=0.1$ & RB & 0.02 & 0.00 & -0.03 & 0.01 & 0.02 & 0.00 & 0.02 &
0.03 \\
& RS & 0.44 & 0.30 & 0.19 & 0.14 & 0.38 & 0.27 & 0.18 & 0.13 \\
& L & 3.6 & 2.7 & 2.8 & 2.3 & 3.5 & 3.0 & 1.7 & 1.9 \\
& U & 2.3 & 3.1 & 2.5 & 2.6 & 3.4 & 3.7 & 1.7 & 2.7 \\
& L${}+{}$U & 5.9 & 5.8 & 5.3 & 4.9 & 6.9 & 6.7 & 3.4 & 4.6 \\[3pt]
$\rho=0.2$ & RB & -0.01 & 0.00 & 0.00 & -0.01 & 0.00 & 0.01 & 0.03 &
0.04 \\
& RS & 0.41 & 0.32 & 0.20 & 0.14 & 0.37 & 0.28 & 0.18 & 0.14 \\
& L & 2.6 & 3.5 & 2.7 & 3.4 & 2.0 & 2.7 & 1.2 & 2.9 \\
& U & 3.0 & 2.5 & 2.7 & 2.5 & 3.6 & 3.5 & 2.8 & 3.2 \\
& L${}+{}$U & 5.6 & 6.0 & 5.4 & 5.9 & 5.6 & 6.2 & 4.0 & 6.1 \\[3pt]
$\rho=0.3$ & RB & -0.01 & 0.01 & 0.02 & 0.01 & -0.01 & 0.00 & 0.00 &
0.03 \\
& RS & 0.43 & 0.32 & 0.20 & 0.15 & 0.38 & 0.28 & 0.18 & 0.14 \\
& L & 3.7 & 3.8 & 2.0 & 3.1 & 2.8 & 2.8 & 2.6 & 2.0 \\
& U & 4.0 & 3.2 & 2.8 & 3.0 & 5.0 & 5.0 & 2.9 & 2.6 \\
& L${}+{}$U & 7.7 & 7.0 & 4.8 & 6.1 & 7.8 & 7.8 & 5.5 & 4.6 \\
\hline
\end{tabular*}
\end{table}

The results obtained for population 3 are presented in Tables~\ref{tabsimres3a} and \ref{tabsimres3b}. We observed no qualitative
difference on populations 1 and 2, and the results for these two
populations are therefore presented in the supplement \cite{cha15} for
brevity. We first consider the results of variance estimation for a
total with the simplified variance estimator $v_{\mathrm{SIMP}}(\hat{Y})$ and
with the bootstrap of PSUs, which are presented in Table~\ref{tabsimres3a}. We note that both variance estimators are
approximately unbiased, with \mbox{absolute} relative biases no greater than
$10 \%$. As expected, $v_{\mathrm{SIMP}}(\hat{Y})$ is slightly negatively biased
while the bootstrap variance estimator is slightly positively biased.
The absolute bias tends to increase with $n_I$, that is, when the
sampling fraction becomes nonnegligible. The simplified variance
estimator is slightly more stable in all scenarios, while the bootstrap
performs slightly better in terms of coverage rates. We now consider
the results obtained for the bootstrap of PSUs when estimating a ratio
and a correlation coefficient, which are presented in Table~\ref{tabsimres3b}. The bootstrap variance estimator is almost unbiased,
with absolute relative biases no greater than $4 \%$. The coverage
rates are well respected in all cases.

\section{Application on the panel for urban policy} \label{applippv}

We illustrate the proposed methods in the context of the Panel for
Urban Policy (PUP), which was conducted by the French General
Secretariat of the Inter-ministerial Committee for Cities (SGCIV). The
PUP is a panel survey in four waves conducted between 2011 and 2014,
which focuses on individuals in the Sensitive Urban Zones (ZUS), and
which collects information on various aspects including security,
employment, precariousness, schooling and health. In this paper, we
focus on the 2011 edition. It~involved two stages of sampling, with the
selection of districts as PSUs, and of households as SSUs. All the
individuals within the selected households were surveyed.

For the purpose of illustration, we consider a subset of districts as
our population $U_I$ of interest. At the first stage, the population
$U_I$ is partitioned into $L=11$ strata $U_{Il}$ according to the
district. In each stratum $U_{Il}$ of size $N_{Il}$, a SI sample
$S_{Il}$ of $n_{Il}$ households is selected and all the individuals
within the households $u_i \in S_{Ih}$ are surveyed. In summary, our
data set consists in a sample of $576$ individuals obtained by
stratified SI cluster sampling of households. The first-stage sampling
rates $f_{Il}=N_{Il}^{-1} n_{Il}$ inside the $L$ strata range from
$0.002$ to $0.017$, which can be considered as negligible.

We are interested in four variables related to health. The variable
$y_1$ gives the perceived health status (very good, good, fair, poor).
The variable $y_2$ is an indicator of chronic disease (with, without).
The variable $y_3$ indicates if the individual is limited by his health
status in his usual activities (very limited, limited, not limited).
The variable $y_4$ indicates if the individual benefits from a free
universal health care (yes, no). For any possible characteristic $c$ of
some variable $y$, we are interested in the proportion
%
\begin{eqnarray}
\label{pc} p_{c} &=& \frac{\sum_{l=1}^L \sum_{u_i \in U_{Il}} Y_{ic}}{\sum_{l=1}^L
\sum_{u_i \in U_{Il}} N_{i}} \qquad \mbox{with }
Y_{ic}=\sum_{k \in
u_i} 1(y_k=c),
\end{eqnarray}
which is estimated by its substitution estimator
%
\begin{eqnarray}
\label{hatpc} && \hat{p}_{c} = \frac{\sum_{l=1}^L N_{Il} n_{Il}^{-1} \sum_{u_i \in
S_{Il}} Y_{ic}}{\hat{N}} \qquad \mbox{with
} \hat{N} \equiv\sum_{l=1}^L
\frac{N_{Il}}{n_{Il}} \sum_{u_i \in S_{Il}} N_{i}.
\end{eqnarray}
For each proportion, we give the normality-based confidence interval.
For that purpose, we adapt the simplified variance estimator in (\ref
{estvarhtsrssimp2}) to the stratified context and make use of the
linearized variable of $p_{c}$. This leads to the variance estimator
%
\begin{eqnarray}
\label{vstwrhatpc} && v_{\mathrm{STWR}}(\hat{p}_{c}) = \sum
_{l=1}^L \frac{N_{Il}^2}{n_{Il}} s_{El}^2
\qquad \mbox{with } s_{El}^2 = \frac{1}{n_{Il}-1} \sum
_{u_i \in S_{Il}} (E_i-\bar{E}_l)^2
\end{eqnarray}
for $\hat{p}_{c}$, with\vspace*{-3pt}
\begin{eqnarray*}
 E_i &=& \frac{1}{\hat{N}} (Y_{ic}-
\hat{p}_{c}) \quad \mbox{and}
\\
\bar {E}_l  &=&
\frac{1}{n_{Il}} \sum_{u_i \in S_{Il}} E_i.
\end{eqnarray*}
For each proportion, we also give the percentile bootstrap and the
Studentized bootstrap confidence intervals, using the with-replacement
bootstrap of PSUs with $D=1000$ resamples. The results with a nominal
one-tailed error rate of $2.5$\% are presented in Table~\ref{resultillust}. The three confidence intervals are very similar in any
case, though the normality-based confidence intervals tend to be
slightly larger.

\begin{table}[t]
\caption{Substitution estimator of the marginal proportions,
normality-based confidence interval (CI), Percentile bootstrap
confidence interval and Studentized bootstrap confidence interval for
four variables} \label{resultillust}
\begin{tabular*}{\tablewidth}{@{\extracolsep{\fill}}lcccc@{}}
\hline
& \multicolumn{4}{c@{}}{\textbf{Perceived health status}} \\[-4pt]
& \multicolumn{4}{l@{}}{\hrulefill}\\
& Very good & Good & Fair & Poor \\
Estimator $\hat{p}_c$ & $0.19$ & $0.43$ & $0.23$ & $0.15$ \\
Normality-based CI & $[0.15,0.24]$ & $[0.38,0.49]$ & $[0.18,0.28]$ &
$[0.10,0.19]$ \\
Percentile bootstrap CI & $[0.15,0.23]$ & $[0.39,0.48]$ & $[0.19,0.27]$ &
$[0.10,0.20]$ \\
Studentized bootstrap CI & $[0.16,0.24]$ & $[0.39,0.48]$ & $[0.19,0.28]$ &
$[0.11,0.21]$ \\[6pt]
& \multicolumn{4}{c@{}}{\textbf{Indicator of chronic disease}} \\[-4pt]
& \multicolumn{4}{l@{}}{\hrulefill}\\
& With & Without & & \\
Estimator $\hat{p}_c$ & $0.28$ & $0.72$ & & \\
Normality-based CI & $[0.23,0.33]$ & $[0.65,0.79]$ & & \\
Percentile bootstrap CI & $[0.24,0.33]$ & $[0.67,0.76]$ & & \\
Studentized bootstrap CI & $[0.24,0.33]$ & $[0.68,0.77]$ & & \\[6pt]
& \multicolumn{4}{c@{}}{\textbf{Limitation in usual activities}} \\[-4pt]
       & \multicolumn{4}{l@{}}{\hrulefill}\\
& Very limited & Limited & Not limited & \\
Estimator $\hat{p}_c$ & $0.09$ & $0.14$ & $0.77$ & \\
Normality-based CI & $[0.05,0.13]$ & $[0.11,0.18]$ & $[0.70,0.84]$ & \\
Percentile bootstrap CI & $[0.06,0.14]$ & $[0.11,0.18]$ & $[0.71,0.81]$ & \\
Studentized bootstrap CI & $[0.06,0.15]$ & $[0.11,0.18]$ & $[0.72,0.82]$ &
\\[6pt]
& \multicolumn{4}{c@{}}{\textbf{Recipient from a free universal health care}} \\[-4pt]
& \multicolumn{4}{l@{}}{\hrulefill}\\
& Yes & No & & \\
Estimator $\hat{p}_c$ & $0.13$ & $0.87$ & & \\
Normality-based CI & $[0.08,0.18]$ & $[0.80,0.94]$ & & \\
Percentile bootstrap CI & $[0.08,0.18]$ & $[0.82,0.92]$ & & \\
Studentized bootstrap CI & $[0.09,0.19]$ & $[0.83,0.92]$ & & \\
\hline
\end{tabular*}
\end{table}

\section{Proofs of results} \label{secproofs}

\subsection{Proof of Theorem~\texorpdfstring{\protect\ref{theotclber}}{3.1}} \label{secprooftclber}

We note $\hat{Y}_{iB} \equiv N_I n_I^{-1} I_i^B \hat{Y}_i$. Under (H2),
we have $\sum_{u_i \in U_i} E|\hat{Y}_{iB}-Y_i|^{2+\delta
}=O(N_I^{2+\delta} n_I^{-1-\delta})$. We obtain
%
\begin{eqnarray}
\label{lyapcond} \frac{\sum_{u_i \in U_i} E|\hat{Y}_{iB}-Y_i|^{2+\delta}}{V(\hat
{Y}_B)^{1+\delta/2}} & = & O\bigl(n_I^{-\delta/2}
\bigr)
\end{eqnarray}
so that the Lyapunov condition is satisfied and (\ref{tclber}) follows
from the central limit theorem for triangular arrays. Noting $\Delta
=N_I^{-2} n_I \{v_B(\hat{Y}_B)-V(\hat{Y}_B)\}$, we have
%
\begin{eqnarray}
\label{espdelt2} E_{\{n_I^B>0\}}\bigl(\Delta^2\bigr) & = &
\frac{1}{1-(1-f_I)^{N_I}} \sum_{k=1}^{n_I}
\operatorname{Pr}\bigl(n_I^B=k\bigr) E_{\{n_I^B=k\}}
\bigl(\Delta^2\bigr),
\end{eqnarray}
where $\operatorname{Pr}(n_I^B=k)=C_{N_I}^k f_I^k (1-f_I)^{N_I-k}$. Using the fact
that conditionally on $n_I^B$, $S_I^B$ may be seen as a simple random
sample of size $n_I^B$ from $U_I$, we have after some algebra that
there exists some constant $C_5$ such that $E_{\{n_I^B=k\}}(\Delta^2)
\leq C_5/k$ for any $k>0$. This leads to
%
\begin{eqnarray}
\label{espdelt22} && E_{\{n_I^B>0\}}\bigl(\Delta^2\bigr) \leq
\frac{C_5}{1-(1-f_I)^{N_I}} \sum_{k=1}^{N_I}
\frac{C_{N_I}^k f_I^k (1-f_I)^{N_I-k}}{k}.
\end{eqnarray}
The term in the right-hand side of (\ref{espdelt22}) tends to $0$
(see Lemma~1.1 in the supplement \cite{cha15}), which leads to (\ref
{estvarconsber0}). To prove (\ref{estvarconsber}), it suffices to
notice that under (H2) there exists some constant $C_6$ such that $E_{\{
n_I^B=0\}}(\Delta^2) = \{N_I^{-2} n_I V(\hat{Y}_B)\}^2 \leq C_6$, and
that $\operatorname{Pr}(n_I^B=0)=(1-f_I)^{N_I}$ tends to $0$.

\subsection{Proof of Theorem \texorpdfstring{\protect\ref{theotclsrs}}{3.2}} \label{prooftheotclsrs}

\begin{lemma} \label{lemcoupling}
Let $X_t$ and $Z_t$ denote two random variables such that
$E(X_t)=E(Z_t)$. Assume that $E(X_t-Z_t)^2=o\{V(X_t)\}$ and that
$V(X_t) \displaystyle\mathop{\longrightarrow}_{t\to\infty} \infty$. Then
%
\begin{eqnarray}
\label{lemcouplingeq1} \bigl\{V(X_t)\bigr\}^{-1}
V(Z_t) & \displaystyle\mathop{\longrightarrow}_{t\to\infty} & 1.
\end{eqnarray}
Also, if for some distribution $\mathcal{L}_0$
%
\begin{eqnarray}
\label{lemcouplingeq2} \bigl\{V(X_t)\bigr\}^{-0.5}\bigl
\{X_t-E(X_t)\bigr\} & \displaystyle\mathop{
\longrightarrow}_{\mathcal{L}} & \mathcal{L}_0,
\end{eqnarray}
then $\{V(Z_t)\}^{-0.5}\{Z_t-E(Z_t)\} \displaystyle\mathop{\longrightarrow}_{\mathcal{L}} \mathcal{L}_0$.
\end{lemma}

The\vspace*{-1.5pt} proof of Lemma~\ref{lemcoupling} is omitted. We take $X_t=\sum_{u_i \in S_I^B} (\hat{Y}_i-\mu_Y)$ and $Z_t=\sum_{u_i \in S_I} (\hat
{Y}_i-\mu_Y)$. Under assumptions (H1) and (H2), Proposition~\ref{prop2}
implies that the assumptions of Lemma~\ref{lemcoupling} are satisfied.
Using the same proof as for (\ref{tclber}) in Theorem~\ref
{theotclber}, it is easily shown that (\ref{lemcouplingeq2}) holds
with $\mathcal{L}_0$ replaced with the standard normal distribution.
This completes the proof.

\subsection{Proof of Theorem \texorpdfstring{\protect\ref{theotclwr}}{4.1}} \label{prooftheotclwr}

Since $\hat{Y}_{\mathrm{WR}}$ is a sum of independent and identically
distributed random variables, (\ref{tclwr}) follows from the classical
central limit theorem for triangular arrays in the i.i.d. case. After
some algebra, we have
\begin{eqnarray*}
V\bigl(s_X^2\bigr) & = & \frac{1}{n_I}
\biggl[E(X_j-\mu_Y)^4-\frac{n_I-3}{n_I-1}
\bigl\{E(X_j-\mu_Y)^2 \bigr\}^2
\biggr].
\end{eqnarray*}
From (H2), there exists $C_{10}$ such that $V(s_X^2) \leq
C_{10} n_I^{-1}$, so that (\ref{estvarconswr}) follows.

\section*{Acknowledgements}
I thank the Commissariat G\'en\'eral \`a l'\'Egalit\'e des Territoires
for supplying the PUP data used in Section~\ref{applippv}. I am
grateful to the Editor, the Associate Editor and the referee for
numerous helpful suggestions which led to an improvement of the paper.
I thank Anne Ruiz-Gazen for helpful comments on earlier versions of
this article.

\begin{supplement}[id=suppA]
\stitle{Supplement to ``Coupling methods for multistage sampling''\\}
\slink[doi,text=10.1214/ 15-AOS1348SUPP]{10.1214/15-AOS1348SUPP} 
\sdatatype{.pdf}
\sfilename{aos1348\_supp.pdf}
\sdescription{The supplement \cite{cha15} contains additional proofs
of Propositions in Section~1, and additional simulation results in
Section~2.}
\end{supplement}


\printaddresses
\end{document}